\documentclass[11pt]{article}
\usepackage{amsfonts}
\usepackage{latexsym,amssymb,amsmath,graphics,cite}
\usepackage{tikz}
\usepackage{fancyhdr}
\usepackage{lipsum}
\usepackage{lineno,slashbox}

\begin{document}
\newcommand{\qed}{\hphantom{.}\hfill $\Box$\medbreak}
\newcommand{\Proof}{\noindent{\bf Proof \ }}
\newtheorem{theorem}{Theorem}[section]
\newtheorem{proposition}[theorem]{Proposition}
\newtheorem{lemma}[theorem]{Lemma}
\newtheorem{corollary}[theorem]{Corollary}
\newtheorem{remark}[theorem]{Remark}
\newtheorem{example}[theorem]{Example}
\newtheorem{definition}[theorem]{Definition}
\newtheorem{construction}[theorem]{Construction}
\renewcommand{\theequation}{\arabic{section}.\arabic{equation}}

\begin{center}
{\Large\bf New $2$-designs from strong difference families\footnote{Supported by NSFC under Grant $11471032$, and Fundamental Research Funds for the Central Universities under Grant $2016$JBM$071$, $2016$JBZ$012$ (T. Feng), Zhejiang Provincial Natural Science Foundation of China under Grant LY17A010008, and Natural Science Foundation of Ningbo under Grant 2016A610079 (X. Wang).}}

\vskip12pt

Simone Costa$^a$, Tao Feng$^b$, Xiaomiao Wang$^c$\\[2ex]
{\footnotesize $^a$Dipartimento DICATAM, Universit\`a degli Studi di Brescia, Via Valotti 9, I-25123 Brescia, Italy}\\
{\footnotesize $^b$Department of Mathematics, Beijing Jiaotong University, Beijing 100044, P. R. China}\\
{\footnotesize $^c$Department of Mathematics, Ningbo University, Ningbo 315211, P. R. China}\\
{\footnotesize simone.costa@unibs.it, tfeng@bjtu.edu.cn, wangxiaomiao@nbu.edu.cn}
\vskip12pt

\end{center}

\vskip12pt

\noindent {\bf Abstract:} Strong difference families are an interesting class of discrete structures which can be used to derive relative difference families. Relative difference families are closely related to $2$-designs, and have applications in constructions for many significant codes, such as optical orthogonal codes and optical orthogonal signature pattern codes. In this paper, with a careful use of cyclotomic conditions attached to strong difference families, we improve the lower bound on the asymptotic existence results of $(\mathbb{F}_{p}\times \mathbb{F}_{q},\mathbb{F}_{p}\times \{0\},k,\lambda)$-DFs for $k\in\{p,p+1\}$. We improve Buratti's existence results for $2$-$(13q,13,\lambda)$ designs and $2$-$(17q,17,\lambda)$ designs, and establish the existence of seven new $2$-$(v,k,\lambda)$ designs for $(v,k,\lambda)\in\{(694,7,2),(1576,8,1),(2025,9,1),(765,9,2),(1845,9,2),(459,9,4)$, $(783,9,4)\}$.



\section{Introduction}

Throughout this paper, sets and multisets will be denoted by curly braces $\{\ \}$ and square brackets $[\ ]$, respectively. When a set or a multiset is fixed with an ordering, we regard it as a sequence, and denote it by $(\ )$. Every union will be understood as {\em multiset union} with multiplicities of elements preserved. $A{\cup}A{\cup}\cdots{\cup}A$ ($h$ times) will be denoted by $\underline{h}A$. If $A$ and $B$ are multisets defined on a multiplicative group, then $A\cdot B$ denotes the multiset $[ab:a\in A,b\in B]$.

For a positive integer $v$, we abbreviate $\{0,1,\dots,v-1\}$ by $\mathbb{Z}_v$ or $I_v$, with the former indicating that a cyclic group of this order is acting.

A {\em $2$-$(v,k,\lambda)$ design} (also called $(v,k,\lambda)$-BIBD or {\em balanced incomplete block design}) is a pair $(V,\cal{A})$ where $V$ is a set of $v$ {\em points} and $\cal A$ is a collection of $k$-subsets of $X$ (called {\em blocks}) such that every
$2$-subset of $X$ is contained in a unique block of $\cal A$.

A powerful idea to obtain $2$-designs is given by the use of difference families and more generally of relative difference families. This idea was implicitly used in many papers (cf. \cite{bb}). The concept of relative difference families was initially put forward by M. Buratti in 1998 \cite{b98}, and systematically developed and discussed by many other authors hereafter (see, for example, \cite{bg,g,m}). Relative difference families play an important role in the construction of optical orthogonal codes (see, for example, \cite{bjlw,fmy,yyl}) and optical orthogonal signature pattern codes \cite{pc}.

Let $(G,+)$ be a group of order $g$ with a subgroup $N$ of order $n$. A $(G,N,k,\lambda)$ {\em relative difference family} (DF), or $(g,n,k,\lambda)$-DF over $G$ relative to $N$, is a family $\mathfrak{B}=[B_1,B_2,\dots,B_r]$ of $k$-subsets of $G$ such that the list
$$\Delta \mathfrak{B}:=\bigcup_{i=1}^r[x-y:x,y\in B_i, x\not=y]=\underline{\lambda}(G\setminus N),$$
i.e., every element of $G\setminus N$ appears exactly $\lambda$ times in the multiset $\Delta \mathfrak{B}$ while it has no element of $N$. The members of $\mathfrak{B}$ are called {\em base blocks} and the number $r$ equals to $\lambda(g-n)/(k(k-1))$. When $G$ is cyclic, we say that the $(g,n,k,\lambda)$-DF is {\em cyclic}. When $N=\{0\}$, a relative difference family is simply called a {\em difference family}.

An {\em automorphism} of a $2$-$(v,k,\lambda)$ design $D=(X,{\cal A})$ is a permutation on
 $X$ leaving $\cal B$ invariant. All automorphisms of $D$ form a group, called
the {\em full automorphism group} of $D$ and denoted by $Aut(D)$. Any subgroup of $Aut(D)$ is called an {\em automorphism group of $D$}. A $2$-$(v,k,\lambda)$ design admitting $\mathbb{Z}_v$ as its automorphism group is called {\em cyclic}. A $2$-$(v,k,\lambda)$ design is said to be {\em $1$-rotational} if it admits an automorphism
consisting of one fixed point and a cycle of length $v-1$. The following proposition reveals the relation between relative difference families and $2$-designs (cf. \cite{b98,b99}).

\begin{proposition}\label{prop:RDF-BIBD-1}
\begin{itemize}
\item[$(1)$] If there exist a $(G,N,k,\lambda)$-DF and a $2$-$(|N|,k,\lambda)$ design, then there exists a $2$-$(|G|,k,\lambda)$ design.
\item[$(2)$] If there exist a $(G,N,k,\lambda)$-DF and a $2$-$(|N|+1,k,\lambda)$ design, then there exists a $2$-$(|G|+1,k,\lambda)$ design.
\end{itemize}
\end{proposition}

\begin{remark}
If the $(G,N,k,\lambda)$-DF is cyclic in Proposition \ref{prop:RDF-BIBD-1}, and the input $2$-design is cyclic in $(1)$ $($or $1$-rotational in $(2)$$)$, then so is the resulting $2$-design.
\end{remark}

The target of this paper is to construct relative difference families via strong difference families and to use them to construct new $2$-designs. Let $\mathfrak{S}=[F_1,F_2,\dots,F_s]$ be a family of multisets of size $k$ of a group $(G,+)$ of order $g$. We say that $\mathfrak{S}$ is a $(G,k,\mu)$ {\em strong difference family}, or a $(g,k,\mu)$-SDF over $G$, if the list
$$\Delta \mathfrak{S}:=\bigcup_{i=1}^s [x-y:x,y\in F_i, x\not=y]=\underline{\mu} G,$$ i.e., every element of $G$ (0 included) appears exactly $\mu$ times in the multiset $\Delta \mathfrak{S}$. The members of $\mathfrak{S}$ are also called {\em base blocks} and the number $s$ equals to $\mu g/(k(k-1))$. Note that $\mu$ is necessarily even since the element $0\in G$ is expressed in even ways as differences in any multiset.

\begin{proposition}\label{prop:nece SDF}
A $(G,k,\mu)$-SDF exists only if $\mu$ is even and $\mu |G|\equiv 0 \pmod{k(k-1)}$.
\end{proposition}

M. Buratti \cite{b99} in 1999 introduced the concept of strong difference families to establish systematic constructions for relative difference families. He named his main construction as ``the fundamental construction" for relative difference families in Theorem 3.1 in \cite{b99}. K. Momihara \cite{m} developed Buratti's technique to give the following theorem. Throughout this paper we always write
\begin{eqnarray}\label{Q(e,m)}
Q(d,m)=\frac{1}{4}(U+\sqrt{U^2+4d^{m-1}m})^2, \mbox{ where } U=\sum_{h=1}^m {m \choose h}(d-1)^h(h-1)
\end{eqnarray}
for given positive integers $d$ and $m$.

\begin{theorem}\label{thm:DF-1}{\rm (Theorems $5.8$ and $5.9$ in \cite{m})}
If there exists a $(G,k,\mu)$-SDF with $\mu=\lambda d$, then there exists a $(G\times \mathbb{F}_q,G\times \{0\},k,\lambda)$-DF
\begin{itemize}
\item for any even $\lambda$ and any prime power $q\equiv 1 \pmod{d}$ with $q>Q(d,k-1)$;
\item for any odd $\lambda$ and any prime power $q\equiv d+1 \pmod{2d}$ with $q>Q(d,k-1)$.
\end{itemize}
\end{theorem}

We remark that M. Buratti and A. Pasotti first presented Theorem \ref{thm:DF-1} for the case of $\lambda=1$ in their Theorem 5.1 in \cite{bp}.

Theorem \ref{thm:DF-1} shows that any $(G,k,\mu)$-SDF can lead to an infinite family of $(G\times \mathbb{F}_q,G\times \{0\},k,\lambda)$-DFs for any admissible sufficiently large prime power $q$. We shall prove that if the initial SDF in Theorem \ref{thm:DF-1} has some particular patterns, then the lower bound on $q$ can be reduced greatly (see Theorems \ref{thm:DF-2}, \ref{thm:DF-4} and \ref{thm:DF-3}). This enables us to obtain new 2-designs with block sizes $13$ and $17$ (see Theorem \ref{thm:13_17}).

On the other hand, we recall that, despite the fact that many authors worked on the existence of $2$-$(v,k,\lambda)$ designs with $6\leq k\leq 9$, there are still many open cases. In this paper we show that, with a careful application of cyclotomic conditions attached to a strong difference family, it is possible to establish existences of $2$-$(v,k,\lambda)$ designs in some open cases. We can establish the existence of $2$-$(v,k,\lambda)$ designs for $(v,k,\lambda)\in\{(694,7,2),(1576,8,1),(2025,9,1),(765,9,2),(1845,9,2),(459,9,4),(783,9,4)\}$.

\section{Seven new 2-designs}

Let $q$ be a prime power. As usual we denote by $\mathbb{F}_q$ the finite field of order $q$ and by $\mathbb{F}^*_q$ its multiplicative group. If $q\equiv 1 \pmod{e}$, then $C_0^{e,q}$ will denote the group of nonzero $e$th powers of $\mathbb{F}_q$ and once a primitive element $\omega$ of $\mathbb{F}_q$ has been fixed, we set $C_i^{e,q}=\omega^i\cdot C_0^{e,q}$ for $i=0,1,\ldots,e-1$. We refer to the cosets $C_0^{e,q},C_1^{e,q},\ldots,C_{e-1}^{e,q}$ of $C_0^{e,q}$ in $\mathbb{F}^*_q$ as the {\em cyclotomic classes of index $e$}.

Let $A$ be a multisubset of $\mathbb{F}^*_q$. If each cyclotomic coset $C_l^{e,q}$ for $l\in I_e$ contains exactly $\lambda$ elements of $A$, then $A$ is said to be a {\em $\lambda$-transversal} for these cosets. If $A$ is a $1$-transversal, $A$ is often referred to as a {\em representative system for the cosets of $C_0^{e,q}$ in $\mathbb{F}_{q}^*$}.

The following lemma is a special case of ``the fundamental construction" for relative difference families in Theorem 3.1 in \cite{b99}. We outline the proof for completeness.

\begin{lemma}\label{lem:DF-1}
Let $q\equiv 1 \pmod{e}$ be a prime power and $d|e$. Let $S$ be a representative system for the cosets of $C_0^{e,q}$ in $C_0^{d,q}$. Write $\mu=\lambda d(q-1)/e$. Suppose that there exists a $(G,k,\mu)$-SDF consisting of base blocks $F_i=(f_{i0},f_{i1},\ldots,$ $f_{i,k-1})$, $1\leq i\leq n$. If there exists a multiset $[\Phi_{1},\Phi_{2},\dots,\Phi_{n}]$ of ordered $k$-subsets of $\mathbb{F}_q^*$ with $\Phi_i=(\phi_{i0},\phi_{i1},\ldots,\phi_{i,k-1})$, $1\leq i\leq n$, such that for each $h\in G$,
$$\bigcup_{i=1}^n[\phi_{ia}-\phi_{ib}: f_{ia}-f_{ib}=h,(a,b)\in I_k\times I_k,a\neq b]=C_0^{e,q}\cdot D_h,$$
where $D_h$ is a $\lambda$-transversal for the cosets of $C_0^{d,q}$ in $\mathbb{F}_q^*$, then $$\mathfrak{F}=[B_{i}\cdot\{(1,s)\}:1\leq i\leq n,s\in S]$$
forms a $(G\times \mathbb{F}_q,G\times \{0\},k,\lambda)$-DF, where
$B_{i}=\{(f_{i0},\phi_{i0}),(f_{i1},\phi_{i1}),
\ldots,(f_{i,k-1},\phi_{i,k-1})\}.$
\end{lemma}

\Proof It is readily checked that
\begin{eqnarray*}
\Delta \mathfrak{F}= \hspace{-3mm}&\hspace{-3mm}& \bigcup_{s\in S}\bigcup_{i=1}^n (\Delta B_{i}\cdot\{(1,s)\})=\bigcup_{s\in S}\bigcup_{i=1}^n [(f_{ia}-f_{ib},(\phi_{ia}-\phi_{ib})\cdot s):(a,b)\in I_k\times I_k,a\neq b]\\
=\hspace{-3mm}&\hspace{-3mm}&\underline{\lambda}(\bigcup_{s\in S}[\{h\}\times(C_0^{e,q}\cdot D_h\cdot\{s\}):h\in G])=\underline{\lambda} (G\times \mathbb{F}_q^*).
\end{eqnarray*} Note that $C_0^{e,q}\cdot D_h\cdot S=\mathbb{F}_q^*$ for each $h\in G$. \qed

\subsection{A $2$-$(694,7,2)$ design and a $2$-$(v,9,4)$ design for $v\in\{459,783\}$}

In this subsection, we shall apply Lemma \ref{lem:DF-1} with $e=q-1$ to present three new 2-designs. When $e=q-1$, $C_0^{e,q}=\{1\}$ and $S=C_0^{d,q}$.

\begin{lemma}\label{lem:SDF(63,7,2)}
There exists a $(\mathbb{Z}_{63},7,2)$-SDF.
\end{lemma}

\Proof Take
\begin{center}
\begin{tabular}{lll}
$F_{1}=[0,4,15,23,37,58,58]$,&
$F_{2}=[0,1,3,7,13,25,39]$,&
$F_{3}=[0,1,3,11,18,34,47]$.
\end{tabular}
\end{center}

\noindent Then $[F_i: 1\leq i\leq 3]$ forms a $(\mathbb{Z}_{63},7,2)$-SDF. \qed

\begin{lemma}\label{lem:eg-DF-1}
There exists a $(\mathbb{Z}_{63}\times \mathbb{F}_{11},\mathbb{Z}_{63}\times \{0\},7,1)$-DF.
\end{lemma}

\Proof Take the $(\mathbb{Z}_{63},7,2)$-SDF from Lemma \ref{lem:SDF(63,7,2)} as the first components of base blocks of the required DF. Let
\begin{center}
\begin{tabular}{lll}
$B_{1}=\{(0,0),(4,3),(15,5),(23,6),(37,8),(58,1),(58,10)\}$,\\
$B_{2}=\{(0,0),(1,2),(3,4),(7,6),(13,1),(25,10),(39,8)\},$\\
$B_{3}=\{(0,0),(1,4),(3,7),(11,9),(18,2),(34,3),(47,5)\}.$\\
\end{tabular}
\end{center}
Then applying Lemma \ref{lem:DF-1} with $k=7$, $d=2$, $q=11$, $e=10$ and $\lambda=1$, we have
$$\mathfrak{F}=[B_{i}\cdot\{(1,s)\}: 1\leq i\leq 3,s\in C_0^{2,11}]$$
 forms a $(\mathbb{Z}_{63}\times \mathbb{F}_{11},\mathbb{Z}_{63}\times \{0\},7,1)$-DF.
It is readily checked that each $D_h$, $h\in \mathbb{Z}_{63}$, is a representative system for the cosets of $C_0^{2,11}$ in $\mathbb{F}_{11}^*$. \qed

\begin{theorem}\label{thm:BIBD(694,7,2)}
There exists a $2$-$(694,7,2)$ design.
\end{theorem}

\Proof By Lemma \ref{lem:eg-DF-1}, there exists a $(693,63,7,1)$-DF, which implies the existence of a $(693,63,7,2)$-DF. Applying Proposition \ref{prop:RDF-BIBD-1}(2) with a 2-$(64,7,2)$ design, which exists by Theorem 2.5 in \cite{a}, we get a 2-$(694,7,2)$ design. \qed

Combining the known result on the existence of $2$-$(v,7,2)$ designs from Theorem 2.5 in \cite{a}, we have the following corollary.

\begin{corollary}\label{cor:BIBD(v,7,2)}
There exists a $2$-$(v,7,2)$ design if and only if $v\equiv 1,7\pmod{21}$ with the definite exception of $v=22$ and possible exceptions of $v\in\{274,358,574,988,994\}$.
\end{corollary}

\begin{lemma}\label{lem:SDF(27,9,8)}
There exists a $(\mathbb{Z}_{27},9,8)$-SDF.
\end{lemma}

\Proof Take
\begin{center}
\begin{tabular}{lll}
$F_1=[0,3,3,8,8,17,17,23,23]$,&
$F_2=[0,1,2,3,19,4,5,8,12]$,\\
$F_3=[0,1,2,3,19,6,11,13,17]$.
\end{tabular}
\end{center}

\noindent Then $[F_i: 1\leq i\leq 3]$ forms a $(\mathbb{Z}_{27},9,8)$-SDF. \qed

\begin{lemma}\label{lem:eg-DF-2}
There exists a $(\mathbb{Z}_{27}\times \mathbb{F}_{q},\mathbb{Z}_{27}\times \{0\},9,4)$-DF for $q\in\{17,29\}$.
\end{lemma}

\Proof  Take the $(\mathbb{Z}_{27},9,8)$-SDF from Lemma \ref{lem:SDF(27,9,8)} as the first components of base blocks of the required DFs.  For $q=17$, let
{\begin{center}
\begin{tabular}{lll}
$B_{1}=\{(0,0),(3,1),(3,16),(8,2),(8,15),(17,3),(17,14),(23,5),(23,12)\}$,\\
$B_{2}=\{(0,0),(1,1),(2,2),(3,7),(19,11),(4,10),(5,5),(8,14),(12,16)\},$\\
$B_{3}=\{(0,0),(1,16),(2,15),(3,10),(19,6),(6,3),(11,2),(13,12),(17,13)\}.$\\
\end{tabular}
\end{center}}

\noindent For $q=29$, let
{\begin{center}
\begin{tabular}{lll}
$B_{1}=\{(0,0),(3,1),(3,28),(8,2),(8,27),(17,3),(17,26),(23,4),(23,25)\}$,\\
$B_{2}=\{(0,0),(1,1),(2,2),(3,4),(19,11),(4,15),(5,5),(8,13),(12,21)\},$\\
$B_{3}=\{(0,0),(1,28),(2,27),(3,25),(19,18),(6,11),(11,19),(13,10),(17,22)\}.$\\
\end{tabular}
\end{center}}

\noindent Then applying Lemma \ref{lem:DF-1} with $k=9$, $d=2$, $e=q-1$ and $\lambda=4$, we have
$$\mathfrak{F}=[B_{i}\cdot\{(1,s)\}: 1\leq i\leq 3,s\in C_0^{2,q}]$$ forms a $(\mathbb{Z}_{27}\times \mathbb{F}_{q},\mathbb{Z}_{27}\times \{0\},9,4)$-DF for $q\in\{17,29\}$.
It is readily checked that each $D_h$, $h\in\mathbb{Z}_{27}$, is a $4$-transversal for the cosets of $C_0^{2,q}$ in $\mathbb{F}_{q}^*$. \qed

\begin{theorem}\label{thm:BIBD(459,9,4)}
There exist a $2$-$(459,9,4)$ design and a $2$-$(783,9,4)$ design.
\end{theorem}

\Proof By Lemma \ref{lem:eg-DF-2}, there exist a $(459,27,9,4)$-DF and a $(783,27,9,4)$-DF. Applying Proposition \ref{prop:RDF-BIBD-1}(1) with a 2-$(27,9,4)$ design, which exists by Theorem 9.1 in \cite{abg}, we get the required 2-designs. \qed

Combining the known result on the existence of $2$-$(v,9,4)$ designs from Theorem 9.1 in \cite{abg},  we have the following corollary.

\begin{corollary}\label{cor:BIBD(v,9,4)}
There exists a $2$-$(v,9,4)$ design if and only if $v\equiv 1,9\pmod{18}$ with the  possible exception of $v=315$.
\end{corollary}

\subsection{A $2$-$(v,9,2)$ design for $v\in\{765,1845\}$}

In this subsection, we shall apply Lemma \ref{lem:DF-1} with $e=(q-1)/2$ to present two new 2-designs. When $e=(q-1)/2$, $C_0^{e,q}=\{1,-1\}$.

\begin{lemma}\label{lem:SDF 45}
There exists a $(\mathbb{Z}_{45},9,8)$-SDF.
\end{lemma}

\Proof Take \begin{center}
\begin{tabular}{lll}
 & $F_1=[0,2,2,15,15,23,23,33,33],$ \\
 &$F_2=F_3=[0,1,4,5,6,7,13,22,33],$ \\
 &$F_4=F_5=[0,2,5,11,21,25,28,36,40].$
\end{tabular}
\end{center}

\noindent Then $[F_i: 1\leq i\leq 5]$ forms a $(\mathbb{Z}_{45},9,8)$-SDF. \qed

\begin{lemma}\label{lem:81_17}
There exists a $(\mathbb{Z}_{45}\times\mathbb{F}_{17},\mathbb{Z}_{45}\times\{0\},9,2)$-DF.
\end{lemma}

\Proof  Take the $(\mathbb{Z}_{45},9,8)$-SDF from Lemma \ref{lem:SDF 45} as the first components of base blocks of the required DF. Let
\begin{center}
\begin{tabular}{l}
$B_1= \{(0,0),(2,1),(2,-1),(15,2),(15,-2),(23,3),(23,-3),(33,5),(33,-5)\},$\\
$B_2= \{(0,0),(1,1),(4,2),(5,3),(6,6),(7,9),(13,4),(22,11),(33,15)\},$\\
$B_4= \{(0,0),(2,3),(5,8),(11,6),(21,12),(25,7),(28,9),(36,2),(40,13)\},$\\
$B_3= B_2\cdot\{(1,-1)\},$\ \ \ \ \ \ \ \ \ \ \ $B_5=B_4\cdot\{(1,-1)\}.$\\
\end{tabular}
\end{center}
Let $S$ be a representative system for the cosets of $C_0^{8,17}=\{1,-1\}$ in $C_0^{2,17}$. Then applying Lemma \ref{lem:DF-1} with $k=9$, $d=2$, $q=17$, $e=8$ and $\lambda=2$, we have
$$\mathfrak{F}=[B_{i}\cdot\{(1,s)\}: 1\leq i\leq 5,s\in S]$$ forms a $(\mathbb{Z}_{45}\times\mathbb{F}_{17},\mathbb{Z}_{45}\times\{0\},9,2)$-DF. It is readily checked that each $D_h$, $h\in\mathbb{Z}_{45}$, is a $2$-transversal for the cosets of $C_0^{2,17}$ in $\mathbb{F}_{17}^*$. \qed

\begin{lemma}\label{lem:81_41}
There exists a $(\mathbb{Z}_{45}\times\mathbb{F}_{41},\mathbb{Z}_{45}\times\{0\},9,1)$-DF.
\end{lemma}

\Proof  Take the $(\mathbb{Z}_{45},9,8)$-SDF from Lemma \ref{lem:SDF 45} as the first components of base blocks of the required DF. Let
\begin{center}
\begin{tabular}{l}
$B_1=\{(0,0),(2,1),(2,-1),(15,2),(15,-2),(23,3),(23,-3),(33,6),(33,-6)\},$\\
$B_2=\{(0,0),(1,1),(4,7),(5,21),(6,12),(7,15),(13,24),(22,4),(33,34)\},$\\
$B_4=\{(0,0),(2,3),(5,31),(11,32),(21,15),(25,9),(28,40),(36,25),(40,35)\},$\\
$B_3= B_2\cdot\{(1,-1)\},$\ \ \ \ \ \ \ \ \ \ \ $B_5=B_4\cdot\{(1,-1)\}.$\\
\end{tabular}
\end{center}
Let $S$ be a representative system for the cosets of $C_0^{20,41}=\{1,-1\}$ in $C_0^{4,41}$. Then applying Lemma \ref{lem:DF-1} with $k=9$, $d=4$, $q=41$, $e=20$ and $\lambda=1$, we have
$$\mathfrak{F}=[B_{i}\cdot\{(1,s)\}: 1\leq i\leq 5,s\in S]$$
forms a $(\mathbb{Z}_{45}\times\mathbb{F}_{41},\mathbb{Z}_{45}\times\{0\},9,1)$-DF. It is readily checked that each $D_h$, $h\in\mathbb{Z}_{45}$, is a representative system for the cosets of $C_0^{4,41}$ in $\mathbb{F}_{41}^*$. \qed

\begin{theorem}\label{thm:BIBD(v,9,2)}
There exist a $2$-$(765,9,2)$ design and a $2$-$(1845,9,2)$ design.
\end{theorem}

\Proof By Lemmas \ref{lem:81_17} and \ref{lem:81_41}, there exist a $(765,45,9,2)$-DF and a $(1845,45,9,2)$-DF (note that two $(1845,45,9,1)$-DFs can produce a $(1845,45,9,2)$-DF). Applying Proposition \ref{prop:RDF-BIBD-1}(1) with a 2-$(45,9,2)$ design, which exists by Theorem 8.2 in \cite{abg}, we get the required 2-designs. \qed

Combining the known result on the existence of $2$-$(v,9,2)$ designs from Theorem 8.2 in \cite{abg}, we have

\begin{corollary}\label{cor:BIBD(v,9,2)}
There exists a $2$-$(v,9,2)$ design if and only if $v\equiv 1,9\pmod{36}$ with the possible exceptions of $v\in\{189,253,505,837,1197,1837\}$.
\end{corollary}

\subsection{A $2$-$(1576,8,1)$ design and a $2$-$(2025,9,1)$ design}

In this subsection, we shall apply Lemma \ref{lem:DF-1} with $e=(q-1)/4$ to present two new 2-designs.

\begin{lemma}\label{lem:SDF(63,8,8)}
There exists a $(\mathbb{Z}_{p},k,8)$-SDF for $(p,k)\in\{(63,8),(81,9)\}$.
\end{lemma}

\Proof For $(p,k)=(63,8)$, take
\begin{center}
\begin{tabular}{lll}
$F_1=[20,20,-20,-20,29,29,-29,-29]$, \\
$F_2=F_3=F_4=F_5=[0,1,3,7,19,34,42,53]$,\\
$F_6=F_7=F_8=F_9=[0,1,4,6,26,36,43,51]$.
\end{tabular}
\end{center}
For $(p,k)=(81,9)$, take
\begin{center}
\begin{tabular}{lll}
 & $F_1=[0,4,4,-4,-4,37,37,-37,-37],$\\
 &$F_2=F_3=F_4=F_5=[0,1,4,6,17,18,38,63,72],$\\
 &$F_6=F_7=F_8=F_9=[0,2,7,27,30,38,53,59,69].$
\end{tabular}
\end{center}
\noindent Then $[F_i: 1\leq i\leq 9]$ forms the required $(\mathbb{Z}_{p},k,8)$-SDF. \qed

\begin{lemma}\label{lem:63 25}
There exists a $(\mathbb{Z}_{p}\times \mathbb{F}_{25},\mathbb{Z}_{p}\times \{0\},k,1)$-DF for $(p,k)\in\{(63,8),(81,9)\}$.
\end{lemma}

\Proof  Take the $(\mathbb{Z}_{p},k,8)$-SDF from Lemma \ref{lem:SDF(63,8,8)} as the first components of base blocks of the required DFs. Take $p(x)=x^2+2x+3$ to be a primitive polynomial of degree $2$ over $\mathbb{F}_{5}$. Then $\mathbb{F}_{25}$ is the splitting field of $p(x)$ over $\mathbb{F}_{5}$ and, denoted by $\omega$ a root of $p(x)=0$, $\omega$ is a generator of the multiplicative group of $\mathbb{F}_{25}$. Let $\xi=\omega^6$. For $(p,k)=(63,8)$, let
\begin{center}
\begin{tabular}{l}
$B_1= \{(20,1),(20,-1),(-20,\xi),(-20,-\xi),(29,\omega),(29,-\omega),(-29,\omega\xi),(-29,-\omega\xi)\},$\\
$B_2= \{(0,0),(1,1),(3,\omega),(7,\omega^2),(19,\omega^3),(34,\omega^4),(42,\omega^{7}),(53,\omega^{10})\},$\\
$B_6= \{(0,0),(1,\omega),(4,\omega^4),(6,\omega^{20}),(26,\omega^{14}),(36,\omega^{12}),(43,\omega^{15}),(51,\omega^{17})\},$\\
$B_3= B_2\cdot\{(1,-1)\},$\ \ \ \ \ \ \ \ \ \ \ $B_4=B_2\cdot\{(1,\xi)\}$,\ \ \ \ \ \ \ \ \ \ \ $B_5= B_2\cdot\{(1,-\xi)\}$,\\
$B_7=B_6\cdot\{(1,-1)\}$,\ \ \ \ \ \ \ \ \ \ \
$B_8= B_6\cdot\{(1,\xi)\},$\ \ \ \ \ \ \ \ \ \ \ $B_9=B_6\cdot\{(1,-\xi)\}$.\\
\end{tabular}
\end{center}
For $(p,k)=(81,9)$, let
\begin{center}
\begin{tabular}{l}
$B_1=\{(0,0),(4,1),(4,-1),(-4,\xi),(-4,-\xi),(37,\omega),(37,-\omega),(-37,\omega\xi),(-37,-\omega\xi)\},$\\
$B_2=\{(0,0),(1,1),(4,\omega),(6,\omega^2),(17,\omega^3),(18,\omega^4),(38,\omega^5),(63,\omega^7),(72,\omega^8)\},$\\
$B_6=\{(0,0),(2,1),(7,\omega^4),(27,\omega^{17}),(30,\omega^{2}),(38,\omega^{18}),(53,\omega^{8}),(59,\omega^{10}),(69,\omega^{14})\},$\\
$B_3=B_2\cdot\{(1,-1)\},$\ \ \ \ \ \ \ \
$B_4= B_2\cdot\{(1,\xi)\},$\ \ \ \ \ \ \ \
$B_5= B_2\cdot\{(1,-\xi)\},$\\
$B_7= B_6\cdot\{(1,-1)\},$\ \ \ \ \ \ \ \
$B_8= B_6\cdot\{(1,\xi)\},$\ \ \ \ \ \ \ \
$B_9=B_6\cdot\{(1,-\xi)\}.$\\
\end{tabular}
\end{center}
Let $S$ be a representative system for the cosets of $C_0^{6,25}=\{1,-1,\xi,-\xi\}$ in $C_0^{2,25}$. Then applying Lemma \ref{lem:DF-1} with $d=2$, $e=6$ and $\lambda=1$, we have
$$\mathfrak{F}=[B_{i}\cdot\{(1,s)\}: 1\leq i\leq 9,s\in S]$$ forms a $(\mathbb{Z}_{p}\times \mathbb{F}_{25},\mathbb{Z}_{p}\times \{0\},k,1)$-DF. It is readily checked that each $D_h$, $h\in\mathbb{Z}_{p}$, is a representative system for the cosets of $C_0^{2,25}$ in $\mathbb{F}_{25}^*$. \qed

\begin{theorem}\label{thm:BIBD(1576)}
There exist a $2$-$(1576,8,1)$ design and a $2$-$(2025,9,1)$ design.
\end{theorem}

\Proof By Lemma \ref{lem:63 25}, there exist a $(1575,63,8,1)$-DF and a $(2025,81,9,1)$-DF. Applying Proposition \ref{prop:RDF-BIBD-1} with a 2-$(64,8,1)$ design and a 2-$(81,9,1)$ design, which exist by Table 3.3 in \cite{ag}, we get the required designs. \qed

According to Table 3.3 in \cite{ag}, Theorem \ref{thm:BIBD(1576)} provides the first examples of a $2$-$(1576,8,1)$ design and a $2$-$(2025,9,1)$ design so far.

\section{DFs from Paley difference multisets}

If a strong difference family only contains one base block, then it is referred to as a {\em difference multiset} (cf. \cite{b99}) or a {\em regular difference cover} (cf. \cite{abms}).

\begin{lemma}\label{lem:SDF-Paley}{\rm \cite{b99}}
\begin{itemize}
\item[$(1)$] Let $p$ be an odd prime power. Then $\{0\}\cup{\underline 2}\mathbb{F}^{\Box}_p$is an $(\mathbb{F}_p,p,p-1)$-SDF $($called Paley difference multiset of the first type$)$.
\item[$(2)$] Let $p\equiv 3\pmod{4}$ be a prime power. Then ${\underline 2}(\{0\}\cup\mathbb{F}^{\Box}_p)$ is an $(\mathbb{F}_p,p+1,p+1)$-SDF $($called Paley difference multiset of the second type$)$.
\end{itemize}
\end{lemma}

\begin{corollary}\label{cor:Paley DS}
Let $p$ be an odd prime power. Let $\lambda$ be any divisor of $p-1$ and $d=(p-1)/\lambda$. Then there exists an $(\mathbb{F}_{p} \times \mathbb{F}_{q},\mathbb{F}_{p}\times \{0\},p,\lambda)$-DF
\begin{itemize}
\item for any even $\lambda$ and any prime power $q\equiv 1 \pmod{d}$ with $q>Q(d,p-1)$;
\item for any odd $\lambda$ and any prime power $q\equiv d+1 \pmod{2d}$ with $q>Q(d,p-1)$.
\end{itemize}
\end{corollary}

\Proof Apply Theorem \ref{thm:DF-1} with the first type Paley $(\mathbb{F}_{p},{p},{p}-1)$-SDF to complete the proof. \qed

\begin{corollary}\label{cor:Paley DS-1}
Let $p\equiv 3\pmod{4}$ be a prime power. Let $\lambda$ be any divisor of $p+1$ and $d=(p+1)/\lambda$. Then there exists an $(\mathbb{F}_{p} \times \mathbb{F}_{q},\mathbb{F}_{p}\times \{0\},p+1,\lambda)$-DF
\begin{itemize}
\item for any even $\lambda$ and any prime power $q\equiv 1 \pmod{d}$ with $q>Q(d,p-1)$;
\item for any odd $\lambda$ and any prime power $q\equiv d+1 \pmod{2d}$ with $q>Q(d,p-1)$.
\end{itemize}
\end{corollary}

\Proof Apply Theorem \ref{thm:DF-1} with the second type Paley $(\mathbb{F}_{p},{p+1},{p}+1)$-SDF to complete the proof. \qed

The lower bounds on $q$ in Corollaries \ref{cor:Paley DS} and \ref{cor:Paley DS-1} are very huge even if $p$ is small. For example by Corollary \ref{cor:Paley DS}, there is an $(\mathbb{F}_{p}\times \mathbb{F}_{q},\mathbb{F}_{p}\times \{0\},p,1)$-DF for any odd prime powers $p$ and $q$ with $q\equiv p \pmod{2(p-1)}$ and $q>Q(p-1,p-1)$. Take $p=13$. Then $Q(12,12)=7.94968\times 10^{27}$. Thus it would be meaningful to develop a new technique to reduce the bound. Actually, Paley difference multisets have good algebraic properties, but the proofs of Corollaries \ref{cor:Paley DS} and \ref{cor:Paley DS-1} do not make use of them.

The main results in the rest part of this paper will rely heavily on the following theorem, which characterizes existences of elements satisfying certain cyclotomic conditions in a finite field.

\begin{theorem}\label{thm:cyclot bound} {\rm \cite{bp,cj}}
Let $q\equiv 1 \pmod{d}$ be a prime power, let $B=\{b_0,b_1,\dots, b_{m-1}\}$ be an arbitrary m-subset of $\mathbb{F}_q$ and let $(\beta_0,\beta_1,\dots,\beta_{m-1})$ be an arbitrary element of $\mathbb{Z}_d^m$. Set $X=\{x\in \mathbb{F}_q: x-b_i\in C^{d,q}_{\beta_i} \mbox{ for } i=0,1,\dots,m-1\}$. Then $X$ is not empty for any prime power $q\equiv 1 \pmod{d}$ and $q>Q(d,m)$.
\end{theorem}

The case of $m=3$ in Theorem \ref{thm:cyclot bound} was first shown by Buratti \cite{b}.
Then a proof similar to that of $m=3$ allows Chang and Ji \cite{cj}, Abel and Buratti \cite{ab0}, and Buratti and Pasotti \cite{bp} to generalize this result to any $m$. Theorem \ref{thm:cyclot bound} is derived from Weil's Theorem (see \cite{ln}, Theorem 5.41) on multiplicative character sums and plays an essential role in the asymptotic existence problem for difference families (cf. \cite{cwz}).

\subsection{The use of Lemma \ref{lem:DF-1} with $e=(q-1)/4$}

The following lemma, which is a corollary of Lemma \ref{lem:DF-1}, gives us a simple but effective starting point to reduce the lower bound on $q$ in Corollary \ref{cor:Paley DS}.

\begin{lemma}\label{lem:DF-2}
Let $p\equiv 1 \pmod{4}$ be a prime power and $\lambda$ be a divisor of $(p-1)/4$. Write $d=(p-1)/4\lambda$. Let $q$ be a prime power satisfying $\lambda(q-1)\equiv 0\pmod{p-1}$. Let $\delta$ be a generator of $C_0^{2,p}$ and $\xi$ be a primitive $4$th root of unity in $\mathbb{F}_{q}$. Take the first type Paley $(\mathbb{F}_{p},p,p-1)$-SDF whose unique base block $(f_{0},f_1,\ldots,f_{p-1})=$
\begin{eqnarray}\label{1st Paley}
(0,\delta,\delta,-\delta,-\delta,\delta^2,\delta^2,-\delta^2,-\delta^2,\ldots,
\delta^{\frac{p-1}{4}},\delta^{\frac{p-1}{4}},-\delta^{\frac{p-1}{4}},-\delta^{\frac{p-1}{4}}).
\end{eqnarray}
Suppose that one can choose an appropriate multiset $(\phi_{0},\phi_1,\ldots,\phi_{p-1})=$
\begin{eqnarray}\label{1st Paley-y}
(0,y_1,-y_1,\xi y_1,-\xi y_1,y_2,-y_2,\xi y_2,-\xi y_2,\ldots,y_{\frac{p-1}{4}},-y_{\frac{p-1}{4}},\xi y_{\frac{p-1}{4}},-\xi y_{\frac{p-1}{4}})\end{eqnarray}
such that $\{y_1,y_2,\ldots,y_{(p-1)/4}\}\subseteq {\mathbb F}_{q}^*$ and for each $h\in \mathbb{F}_{p}$,
\begin{eqnarray*}\label{1st Paley-dh}
[\phi_{a}-\phi_{b}:f_{a}-f_{b}=h,(a,b)\in I_{p}\times I_{p},a\neq b]=\{1,-1,\xi,-\xi\}\cdot D_h,
\end{eqnarray*}
where $D_h$ is a $\lambda$-transversal for the cosets of $C_0^{d,q}$ in $\mathbb{F}_{q}^*$. Let $S$ be a representative system for the cosets of $\{1,-1,\xi,-\xi\}$ in $C_0^{d,q}$. Let
$
B=\{(f_0,\phi_0),(f_1,\phi_1),\ldots,(f_{p-1},\phi_{p-1})\}.
$
Then
$$\mathfrak{F}=[B\cdot\{(1,s)\}:s\in S]$$
forms an $(\mathbb{F}_{p}\times \mathbb{F}_{q},\mathbb{F}_{p}\times \{0\},p,\lambda)$-DF. \end{lemma}

\Proof For $p\equiv 1 \pmod{4}$, $\delta^{(p-1)/4}=-1$ in $\mathbb{F}_{p}$, so the first type Paley $(\mathbb{F}_{p},p,p-1)$-SDF can be written as (\ref{1st Paley}). By the definition of $d$, $d$ and $4$ are both divisors of $(p-1)/\lambda$, so combining the assumption $\lambda(q-1)\equiv 0\pmod{p-1}$, we have $d$ and $4$ are both divisors of $q-1$. This makes $C_0^{d,q}$ and $\xi$ meaningful.
Since $q$ is odd, $y_i\neq -y_i$ and $\xi y_i\neq -\xi y_i$ for any $1\leq i\leq (p-1)/4$, $B\cdot\{(1,s)\}$ is a set of size $p$ for any $s\in S$. Then apply Lemma \ref{lem:DF-1} with $n=1$, $k=p$, $d=(p-1)/4\lambda$ and $e=(q-1)/4$ to obtain the required $(\mathbb{F}_{p}\times \mathbb{F}_{q},\mathbb{F}_{p}\times \{0\},p,\lambda)$-DF. Note that $\{1,-1,\xi,-\xi\}=C_0^{(q-1)/4,q}$. \qed

In order to facilitate the reader to understand the use of Lemma \ref{lem:DF-2}, we shall begin with the first type Paley $(\mathbb{F}_{p},p,p-1)$-SDFs for $p\in\{13,17\}$.

\begin{theorem}\label{eg:13_1}
There exists an $(\mathbb{F}_{13}\times\mathbb{F}_q,\mathbb{F}_{13}\times\{0\},13,1)$-DF for all primes $q\equiv 1 \pmod{12}$ with the possible exceptions of
$q\in E_{13}=\{37, 61, 73, 97, 109, 181, 313, 337, 349,$ $ 373, 409, 421, 541, 577, 829, 853, 1129, 1741, 2473\}.$
\end{theorem}

\Proof To apply Lemma \ref{lem:DF-2}, take the first type Paley $(\mathbb{F}_{13},13,12)$-SDF $(0,4,4,-4,-4,3$, $3,-3,-3,-1,-1,1,1)$ as the first components of base blocks of the required DF. Let $\xi$ be a primitive $4$th root of unity in $\mathbb{F}_q$ and
\begin{eqnarray*}
B=\{ \hspace{-3.5mm}&\hspace{-3.5mm}&(0,0),(4,y_1),(4,-y_1),(-4,y_1\xi),(-4,-y_1\xi),(3,y_2),(3,-y_2),\\
\hspace{-3.5mm}&\hspace{-3.5mm}&(-3,y_2\xi),(-3,-y_2\xi),(-1,y_3),(-1,-y_3),(1,y_3\xi),(1,-y_3\xi)\}.
\end{eqnarray*}
Since $D_h=D_{13-h}$, $h\in\mathbb{F}_{13}$, we only need to consider the differences $D_h$ for $0\leq h\leq ¡Ü6$. It is readily checked that w.l.o.g.,
\begin{center}
\begin{tabular}{ll}
$D_0=2\cdot [y_1,y_2,y_3],$&
$D_1=[y_2-y_1,y_2+y_1,y_3],$\\
$D_2=[y_3(1-\xi),y_3-y_2\xi,y_3+y_2\xi],$&
$D_3=[y_2,y_3-y_1\xi,y_3+y_1\xi],$\\
$D_4=[y_1,y_3-y_2,y_3+y_2],$ &
$D_5=[y_1(1-\xi),y_3-y_1,y_3+y_1],$\\
$D_6=[y_2(1-\xi),y_2-y_1\xi,y_2+y_1\xi].$
\end{tabular}
\end{center}
Now we need to pick up appropriate $y_1, y_2, y_3$ in $\mathbb{F}_q^*$ such that each $D_h$, $h\in\mathbb{F}_{13}$, is a representative system for the cosets of $C_0^{3,q}$ in $\mathbb{F}_q^*$.

We here list one possible scheme to pick up $y_1, y_2, y_3$ in the following table according to cyclotomic classes of index $3$ which $1-\xi$ belongs to (for instance, if $1-\xi\in C_0^{3,q}$, then the second column of this table requires $y_1\in C_0^{3,q}$, $y_2\in C_1^{3,q}$, $y_2-y_1\in C_0^{3,q}$, and so on). One can check that if $y_1, y_2$ and $y_3$ satisfy the cyclotomic conditions in the following table, then each $D_h$, $h\in\mathbb{F}_{13}$, is a representative system for the cosets of $C_0^{3,q}$ in $\mathbb{F}_q^*$.

\begin{table}[h]
\begin{center}
\begin{tabular}{cccccc|c|c|c|c|}

\cline{1-4}\cline{7-10} \multicolumn{1}{|c|}{$1-\xi$} &\multicolumn{1}{|c|}{$0$}& \multicolumn{1}{|c|}{$1$}& \multicolumn{1}{|c|}{$2$}&&&$1-\xi$ &$ 0$& $1$& $2$ \\
\cline{1-4}\cline{7-10} \multicolumn{1}{|c|}{$y_1$} & \multicolumn{1}{|c|}{$0$} &\multicolumn{1}{|c|}{$0$} & \multicolumn{1}{|c|}{$0$} &&& $y_3$ & $2$ &$2 $ & $2$\\
\cline{1-4} \multicolumn{1}{|c|}{$y_2$} & \multicolumn{1}{|c|}{$1 $} &\multicolumn{1}{|c|}{$1 $} & \multicolumn{1}{|c|}{$1 $} &&& $y_3-y_1$ & $1$&$0 $ & $0 $\\
\multicolumn{1}{|c|}{$y_2-y_1$} & \multicolumn{1}{|c|}{$0$} &\multicolumn{1}{|c|}{$0$} & \multicolumn{1}{|c|}{$0$}  &&&$y_3+y_1$ & $2$&$2 $ & \multicolumn{1}{|c|}{$1 $}\\
\multicolumn{1}{|c|}{$y_2+y_1$} & \multicolumn{1}{|c|}{$1 $}&\multicolumn{1}{|c|}{$1 $} & \multicolumn{1}{|c|}{$1 $} &&&$y_3-y_2$ & $1$ &$1 $ & $1 $\\
\multicolumn{1}{|c|}{$y_2-y_1\xi$} & \multicolumn{1}{|c|}{$0$}&\multicolumn{1}{|c|}{$0$} & \multicolumn{1}{|c|}{$1 $} &&&$y_3+y_2$ & $2$ &$2 $ & $2 $\\
\multicolumn{1}{|c|}{$y_2+y_1\xi$} & $2$& \multicolumn{1}{|c|}{$1 $} & \multicolumn{1}{|c|}{$2 $} &&&$y_3-y_1\xi$ & $0$ &$0 $ & $0 $\\
\cline{1-4} &&&&&& $y_3+y_1\xi$ & $2$&$2 $ & \multicolumn{1}{|c|}{$2 $}\\
 &&&&&& $y_3-y_2\xi$ & $0$ &$1 $ & $0 $\\
&&&&&& $y_3+y_2\xi$ & $1$ &$2 $ & $2 $\\
\cline{7-10}
\end{tabular}
\end{center}
\end{table}

Next it suffices to apply Theorem \ref{thm:cyclot bound} to pick up the required $y_1, y_2, y_3$ step by step. First pick up any element of $C_0^{3,q}$ as $y_1$. Once $y_1$ is fixed, by the above table (Row 3 to Row 7 in the left side) there are five cyclotomic conditions on $y_2$, four of which are related with $y_1$, so by Theorem \ref{thm:cyclot bound}, such $y_2$ exists for all primes $q\equiv 1 \pmod{12}$ and $q>Q(3,5)=323433$. Similarly, once $y_1$ and $y_2$ are fixed, there are nine cyclotomic conditions on $y_3$, so by Theorem \ref{thm:cyclot bound} again, such $y_3$ exists for all primes $q\equiv 1 \pmod{12}$ and $q>Q(3,9)=9.68583\times 10^9$.

On the other hand we have checked by computer search that there exist $y_1, y_2, y_3$ in $\mathbb{F}_q^*$ for all primes $q\equiv 1 \pmod{12}$, $q\leq Q(3,9)$ and $q\not\in E_{13}$ such that each $D_h$, $h\in\mathbb{F}_{13}$, is a representative system for the cosets of $C_0^{3,q}$ in $\mathbb{F}_q^*$ (note that it is not necessary to require $y_1, y_2, y_3$ to always satisfy the above table). The interested reader may get a copy of these data from the authors. \qed

\begin{theorem}\label{eg:17_1}
There exists an $(\mathbb{F}_{17}\times\mathbb{F}_q,\mathbb{F}_{17}\times\{0\},17,1)$-DF for all primes $q\equiv 1 \pmod{16}$ and $q>Q(4,13)=3.44807\times 10^{17}$, or $q\in S_{17}\cup\{p: p\ {\rm is\ a\ prime}, p\equiv1\pmod{16}, 6673\leq p\leq 9857\}$, where $S_{17}=$
\begin{center}
\begin{tabular}{ll}
&$\{17, 881, 1297, 1601, 1873, 2017, 2129, 2657, 2753, 2801, 2897, 3089, 3121, 3217$, \\&
$3313, 3361, 3617, 3697, 3761, 3793, 3889, 4001, 4049, 4129, 4241, 4273, 4289, 4481, $\\&
$4561, 4657, 4721, 4801, 4817, 4993, 5009, 5233, 5281, 5297, 5393, 5441, 5521, 5569, $\\&
$5857, 5953, 6113,  6257, 6337, 6449, 6529 \}$.
\end{tabular}
\end{center}
\end{theorem}

\Proof To apply Lemma \ref{lem:DF-2}, take the first type Paley $(\mathbb{F}_{17},17,16)$-SDF: $(0,-8,-8,8,8,$ $-4,-4,4,4,-2,-2,2,2,-1,-1,1,1)$. Let $\xi$ be a primitive $4$th root of unity in $\mathbb{F}_q$ and
\begin{eqnarray*}
B=\{ \hspace{-3.5mm}&\hspace{-3.5mm}&(0,0),(-8,y_1),(-8,-y_1),(8,y_1\xi),(8,-y_1\xi),(-4,y_2),(-4,-y_2),(4,y_2\xi),(4,-y_2\xi),\\
\hspace{-3.5mm}&\hspace{-3.5mm}&(-2,y_3),(-2,-y_3),(2,y_3\xi),(2,-y_3\xi),(-1,y_4),(-1,-y_4),(1,y_4\xi),(1,-y_4\xi)\}.
\end{eqnarray*}
Since $D_h=D_{17-h}$, $h\in\mathbb{F}_{17}$, we only need to consider the differences $D_h$ for $0\leq h\leq ¡Ü8$. It is readily checked that w.l.o.g.,
\begin{center}
\begin{tabular}{ll}
$D_0=2\cdot [y_1,y_2,y_3,y_4],$&
$D_1=[y_1(1-\xi),y_4,y_4-y_3,y_4+y_3],$\\
$D_2=[y_3,y_3-y_2,y_3+y_2,y_4(1-\xi)],$&
$D_3=[y_4-y_2,y_4+y_2,y_4-y_3\xi,y_4+y_3\xi],$\\
$D_4=[y_2,y_2-y_1,y_2+y_1,y_3(1-\xi)],$ &
$D_5=[y_2-y_1\xi,y_2+y_1\xi,y_4-y_2\xi,y_4+y_2\xi],$\\
$D_6=[y_3-y_1,y_3+y_1,y_3-y_2\xi,y_3+y_2\xi],$ &
$D_7=[y_3-y_1\xi,y_3+y_1\xi,y_4-y_1,y_4+y_1],$\\
$D_8=[y_1,y_2(1-\xi),y_4-y_1\xi,y_4+y_1\xi].$
\end{tabular}
\end{center}
Now we need to pick up appropriate $y_1, y_2, y_3, y_4$ in $\mathbb{F}_q^*$ such that each $D_h$, $h\in\mathbb{F}_{17}$, is a representative system for the cosets of $C_0^{4,q}$ in $\mathbb{F}_q^*$.
We here list one possible scheme to pick up $y_1, y_2, y_3, y_4$ in the following table according to cyclotomic classes of index $4$ which $1-\xi$ belongs to. Similar arguments to those in Theorem \ref{eg:13_1} show that for all primes $q\equiv 1 \pmod{16}$ and $q>Q(4,13)=3.44807\times 10^{17}$, these required $y_1, y_2, y_3, y_4$ exist by iterated application of Theorem \ref{thm:cyclot bound}.

\begin{table}[h]
\begin{center}
\begin{tabular}{|c|c|c|c|c|cc|c|c|c|c|c|}
\cline{1-5}\cline{8-12} $1-\xi$ &$ 0$& $1$& $2$ & $3$&&&$1-\xi$ &$ 0$& $1$& $2$ & $3$\\
\cline{1-5}\cline{8-12} {$y_1$} & $ 0$ &$0 $ & $0 $ & $1$&&&$y_4$ & $3$ &$3 $ & $3$ &$2$\\
\cline{1-5}  {$y_2$} & $1$ &$1 $ & $1 $ & $0 $&&&$y_4-y_1$ & $2$&$2 $ & $2 $& $2$\\
{$y_2-y_1$} & $0$ &$0 $ & $2$  & $1$&&&$y_4+y_1$ & $3$&$3 $ & $3 $& $3$\\
{$y_2+y_1$} & $3$&$2 $ & $3 $ & $3$&&&$y_4-y_2$ & $0$ &$0 $ & $0 $& $0$\\
{$y_2-y_1\xi$} & $0$&$0 $ & $0 $ & $0$&&&$y_4+y_2$ & $1$ &$1 $ & $1 $& $1$\\
$y_2+y_1\xi$ & $1$& $1 $ & $1 $ & $1$&&&$y_4-y_3$ & $1$ &$0 $ & $0 $& $1$\\
\cline{1-5} $y_3$ & $2$ &$2 $ & $2$ &$3$&&&$y_4+y_3$ & $2$ &$2 $ & $1 $& $3$\\
$y_3-y_1$ & $0$&$0 $ & $0 $& $0$&&&$y_4-y_1\xi$ & $2$ &$1 $ & $1 $& $0$\\
$y_3+y_1$ & $1$&$1 $ & $1 $& $1$&&&$y_4+y_1\xi$ & $3$&$3 $ & $2 $& $2$\\
$y_3-y_2$ & $0$ &$1 $ & $0 $& $0$&&&$y_4-y_2\xi$ & $2$ &$2 $ & $2 $ & $2$\\
$y_3+y_2$ & $1$ &$3 $ & $3 $& $2$&&&$y_4+y_2\xi$ & $3$ &$3 $ & $3 $& $3$\\
$y_3-y_1\xi$ & $0$ &$0 $ & $0 $& $0$&&&$y_4-y_3\xi$ & $2$ &$2 $ & $2 $ & $2$\\
$y_3+y_1\xi$ & $1$&$1 $ & $1 $& $1$&&&$y_4+y_3\xi$ & $3$ &$3 $ & $3 $& $3$\\
\cline{8-12}$y_3-y_2\xi$ & $2$ &$2 $ & $2 $& $2$ \\
$y_3+y_2\xi$ & $3$ &$3 $ & $3 $& $3$\\
\cline{1-5}
\end{tabular}
\end{center}
\end{table}

On the other hand, with the aid of computer, we can pick up the required $y_1, y_2, y_3, y_4$ in $\mathbb{F}_q^*$ for $q\in S_{17}\cup\{p: p\ {\rm is\ a\ prime}, p\equiv1\pmod{16}, 6673\leq p\leq 9857\}$. The interested reader may get a copy of these data from the authors. \qed

Emulating the proofs of Theorems \ref{eg:13_1} and \ref{eg:17_1}, we shall establish asymptotic existences of $(\mathbb{F}_{p}\times \mathbb{F}_{q},\mathbb{F}_{p}\times \{0\},p,\lambda)$-DFs for all prime powers $p\equiv 1\pmod{4}$ in Theorem \ref{thm:DF-2}. First we need to make a general analysis of $D_h$, $h\in \mathbb{F}_{p}$. Lemmas \ref{lem:DF-2-Dh} and \ref{lem:DF-2-Dh1} follow the notation in Lemma \ref{lem:DF-2}.

\begin{lemma}\label{lem:DF-2-Dh}
For each $h\in \mathbb{F}_{p}^*$, let
$$T_h=[\phi_{a}-\phi_{b}:f_{a}-f_{b}=h,(a,b)\in I_{p}\times I_{p},a\neq b].$$
Then $T_h=\{1,-1,\xi,-\xi\}\cdot D_h$ for some $D_h\subset \mathbb{F}_{q}$ and the size of $D_h$ is $(p-1)/4$. Furthermore, $D_h=D_{-h}$ and w.l.o.g., $D_h$ consists of elements having the following types:
\begin{center}
\begin{tabular}{llll}
$(I)$ $y_i-y_j$, $y_i+y_j$; & $(II)$ $y_i-y_j\xi$, $y_i+y_j\xi$; &
$(III)$ $y_i(1-\xi)$; & $(IV)$ $y_i$,
\end{tabular}
\end{center}
for some $y_i$ and $y_j$ from $(\ref{1st Paley-y})$ in Lemma $\ref{lem:DF-2}$. \end{lemma}

\Proof The base block of the first type Paley $(\mathbb{F}_{p},p,p-1)$-SDF is of the form (\ref{1st Paley}), so its each nonzero difference $h$ is of the form $\pm(\delta^i-\delta^j)$, $\pm(\delta^i+\delta^j)$, $\pm 2\delta^i$ or $\pm \delta^i$ for some $1\leq i, j\leq (p-1)/4$.

If $h=\pm(\delta^i-\delta^j)$, then $T_h\supset \{1,-1,\xi,-\xi\}\cdot\{y_i-y_j,y_i+y_j\}$.

If $h=\pm(\delta^i+\delta^j)$, then $T_h\supset \{1,-1,\xi,-\xi\}\cdot\{y_i-y_j\xi,y_i+y_j\xi\}$.

If $h=\pm 2\delta^i$, then $T_h\supset \{1,-1,\xi,-\xi\}\cdot\{y_i(1-\xi)\}$.

If $h=\pm \delta^i$, then $T_h\supset \{1,-1,\xi,-\xi\}\cdot\{y_i\}$.

Thus $T_h=\{1,-1,\xi,-\xi\}\cdot D_h$ for some $D_h\subset \mathbb{F}_{q}$, $D_h=D_{-h}$ and $D_h$ consists of elements having Types (I)-(IV). Since each $h\in \mathbb{F}_{p}^*$ occurs $p-1$ times as differences of the Paley SDF, the size of $D_h$ is $(p-1)/4$. \qed

\begin{lemma}\label{lem:DF-2-Dh1}
 \begin{itemize}
\item[$(1)$] W.l.o.g., $D_0=2\cdot\{y_1,y_2,\dots,y_{(p-1)/4}\}$.
\item[$(2)$] Let $p\equiv 5\pmod{8}$.
\begin{itemize}
\item If $h\in C_{0}^{2,p}$, then $D_h$ consists of exactly one element of Type $(IV)$ and $(p-5)/4$ elements of Types $(I)$ and $(II)$.
\item If $h\in C_{1}^{2,p}$, then $D_h$ consists of exactly one element of Type $(III)$ and $(p-5)/4$ elements of Types $(I)$ and $(II)$.
\end{itemize}
\item[$(3)$] Let $p\equiv 1\pmod{8}$.
\begin{itemize}
\item If $h\in C_{0}^{2,p}$, then $D_h$ consists of exactly one element of Type $(III)$, exactly one element of Type $(IV)$ and $(p-9)/4$ elements of Types $(I)$ and $(II)$.
\item If $h\in C_{1}^{2,p}$, then $D_h$ consists of $(p-1)/4$ elements of Types $(I)$ and $(II)$.
\end{itemize}
\item[$(4)$] Let $p\equiv 1\pmod{8}$ and $3\nmid p$. If $h\in C_{0}^{2,p}$, then $D_h$ does not contain elements of the form $y_i$ and $y_i(1-\xi)$ at the same time.
\item[$(5)$] Let $p\equiv 1\pmod{8}$ and $3\mid p$. For each $h\in C_{0}^{2,p}$, $D_h$ must contain elements of the form $y_i$ and $y_i(1-\xi)$ at the same time.
\item[$(6)$] Let $p\equiv 1\pmod{4}$ and $T$ be a representative system for the cosets of $\{1,-1\}$ in $\mathbb{F}_{p}^*$. Any element of Types $(I)$ and $(II)$ must be contained in a unique $D_h$ for some $h\in T$ $($note that the term ``element'' here is a symbolic expression; for example $y_1-y_2$ and $y_3+y_4$ are different element but they may have the same value$)$.
\end{itemize}
\end{lemma}

\Proof The verifications for $(1)$, $(5)$ and $(6)$ are straightforward.

$(2)$ For $p\equiv 5\pmod{8}$, $(p-1)/4$ is odd. Since the size of $D_h$ is $(p-1)/4$, each $D_h$, $h\in {\mathbb F}_{p}^*$, must contain at least one element of Type (III) or (IV).

Examining the sequences (\ref{1st Paley}) and (\ref{1st Paley-y}), for any $h=\pm\delta^i\in C_{0}^{2,p}$, $1\leq i\leq (p-1)/4$, we have $y_{i}\in D_h$. Thus each element of Type (IV) is in $D_h$ for some $h\in C_{0}^{2,p}$, and any $D_h$ for $h\in C_{0}^{2,p}$ contains only one element of Type (IV).

On the other hand, for any $h=\pm2\delta^i\in C_{1}^{2,p}$ (note that $2\in C_{1}^{2,p}$ for $p\equiv 5\pmod{8}$), $1\leq i\leq (p-1)/4$, we have $y_{i}(1-\xi)\in D_h$. Thus each element of Type (III) is in $D_h$ for some $h\in C_{1}^{2,p}$, and any $D_h$ for $h\in C_{1}^{2,p}$ contains only one element of Type (III).

$(3)$ A similar argument to that in (2) to complete the proof (note that $2\in C_{0}^{2,p}$ for $p\equiv 1\pmod{8}$).

$(4)$ Any element of the form $y_i$ is in $D_h$ for some $h=\pm\delta^i$, and any element of the form $y_i(1-\xi)$ is in $D_{h'}$ for some $h'=\pm2\delta^i$. If $h=h'$, then $\delta^i=\pm 2\delta^i$, which is impossible since $p$ is not a power of $3$. \qed

\begin{theorem}\label{thm:DF-2}
Let $\lambda$ be any divisor of $(p-1)/4$ and $d=(p-1)/4\lambda$.
\begin{itemize}
\item[$(1)$] There exists an $(\mathbb{F}_{p}\times \mathbb{F}_{q},\mathbb{F}_{p}\times \{0\},p,\lambda)$-DF for any prime powers $p$ and $q$ with $p\equiv 1,5 \pmod{12}$, $\lambda(q-1)\equiv 0\pmod{p-1}$ and $q>Q(d,p-4)$.
\item[$(2)$] Let $p$ be a power of $9$, and $\xi$ be a primitive $4$th root of unity in $\mathbb{F}_{q}$. If $\lambda>1$ or $1-\xi\not\in C_{0}^{d,q}$, then there exists an $(\mathbb{F}_{p}\times \mathbb{F}_{q},\mathbb{F}_{p}\times \{0\},p,\lambda)$-DF for any prime power $q$ with $\lambda(q-1)\equiv 0\pmod{p-1}$ and $q>Q(d,p-4)$.
\end{itemize}
\end{theorem}

\Proof When $p=5$, $\lambda=d=1$. Since $C_0^{1,q}$ is just $\mathbb{F}_{q}^*$, each $D_h$ contains only one element. Then the conclusion is straightforward by Lemma \ref{lem:DF-2}.

When $p\equiv 5\pmod{8}$ is a prime power and $p>5$, by Lemma \ref{lem:DF-2-Dh1}(1), (2) and (6), applying Theorem \ref{thm:cyclot bound}, one can always pick up appropriate $y_1,y_2,\ldots,y_{(p-1)/4}$ such that each $D_h$, $h\in \mathbb{F}_p$, is a $\lambda$-transversal for the cosets of $C_0^{d,q}$ in $\mathbb{F}_{q}^*$ for any divisor $\lambda$ of $(p-1)/4$ and any prime power $q$ with $\lambda(q-1)\equiv 0\pmod{p-1}$ and $q>Q(d,p-4)$. Then apply Lemma \ref{lem:DF-2} to complete the proof.

When $p\equiv 1\pmod{8}$ is a prime power, by Lemma \ref{lem:DF-2-Dh1}(3) and (6), to apply Theorem \ref{thm:cyclot bound}, one thing to be careful of is to examine $D_h$'s, $h\in C_{0}^{2,p}$, of the form $[y_{i_h},y_{j_h}(1-\xi),\cdots]$ such that each of these $D_h$'s is a $\lambda$-transversal for the cosets of $C_0^{d,q}$ in $\mathbb{F}_{q}^*$. When $\lambda>1$, it is easy to do it. It suffices to consider the case of $\lambda=1$.

If $i_{h}=j_{h}$ for some $h\in C_{0}^{2,p}$, then by Lemma \ref{lem:DF-2-Dh1}(4), $p$ must be a power of $9$, and then by Lemma \ref{lem:DF-2-Dh1}(5), for each $h\in C_{0}^{2,p}$, $i_{h}=j_{h}$. In this case, if $1-\xi\not\in C_{0}^{d,q}$, then $y_{i_{h}}$ and $y_{i_{h}}(1-\xi)$ are not in the same coset. Thus we can apply Theorem \ref{thm:cyclot bound} and Lemma \ref{lem:DF-2} to complete the proof.

If $i_{h}\neq j_{h}$ for any $h\in C_{0}^{2,p}$, then $y_{i_h}\neq y_{j_h}$ for any $h\in C_{0}^{2,p}$. By Lemma \ref{lem:DF-2-Dh1}(5), $p\neq 9$, so $p\geq 17$. Let $1-\xi\in C_{a}^{d,q}$ for some $0\leq a<d$. Note that $\{y_{i_{h}}:h\in\{\delta,\delta^2,\ldots,\delta^{d}\}\}=
\{y_{j_{h}}:h\in\{\delta,\delta^2,\ldots,\delta^{d}\}\}=\{y_1,y_2,\ldots,y_d\}$.
To apply Theorem \ref{thm:cyclot bound}, we need to find a bijection $\pi:\{y_1,y_2,\dots,y_{d}\}\rightarrow \mathbb{Z}_{d}$ such that $\pi(y_{i_{h}})\not\equiv \pi(y_{j_{h}})+a\pmod{d}$ for any $D_h=[y_{i_h},y_{j_h}(1-\xi),\cdots]$, $h\in\{\delta,\delta^2,\ldots,\delta^{d}\}$.

Let $\alpha$ be the permutation of $\mathbb{Z}_d$ such that $\alpha(y_{i_h})=y_{j_h}$ where $y_{i_h},y_{j_h}(1-\xi)\in D_h$ and $i_{h}\neq j_{h}$. Thus $h=\delta^{i_h}=\pm 2\delta^{j_h}$. Clearly $\alpha$ is a product of cycles of the same length $l$ that is the multiplicative order of $2$ in $\mathbb{F}_{p}^*/\{1,-1\}$ and $l\geq 2$. Let us write $\alpha$ as $(y_{1,0},y_{1,1},\dots, y_{1,l-1})(y_{2,0},y_{2,1},\dots, y_{2,l-1})\cdots (y_{u,0},y_{u,1},\dots, y_{u,l-1})$ where $u=d/l$.
If $l\geq 3$, or $l=2$ and $u\nmid a$, then a simple number theory argument shows that there exists an element $r\in \mathbb{Z}_d$ such that $ru\not\equiv -a\pmod{d}$ and $(r,l)=1$. We set $\pi(y_{i,j})=i+ruj\pmod{d}$ for all $1\leq i\leq u$ and $0\leq j\leq l-1$. Then $\pi$ is a bijection from the $y_{i,j}$'s onto $\mathbb{Z}_d$ such that $\pi(y_{i,j})-\pi(y_{i,(j+1)})\not\equiv a\pmod{d}$, where the arithmetic $j+1$
is reduced modulo $l$. If $l=2$ and $u\mid a$, then $u=d/2$, which implies $a=0$ or $d/2$. Since $p\geq 17$, $d=(p-1)/4\geq 4$. So $a\neq 1$ and $a\neq d-1$. We set $\pi(y_{i,0})=2i$ and $\pi(y_{i,1})=2i-1$ for all $1\leq i\leq u$. Then $\pi$ is a bijection from the $y_{i,0}$'s and $y_{i,1}$'s onto $\mathbb{Z}_d$ such that $\pi(y_{i,0})-\pi(y_{i,1})\not\equiv \pm a\pmod{d}$. Thus we can apply Theorem \ref{thm:cyclot bound} and Lemma \ref{lem:DF-2} to complete the proof. \qed

Compared with Corollary \ref{cor:Paley DS}, Theorem \ref{thm:DF-2} not only reduces the lower bound on $q$ but also relax the congruence condition on $q$ in some circumstances. Actually, it is easy to see that $Q(e,m)<Q(e,m')$ for $m<m'$ and $Q(e,m)<Q(e',m)$ for $e<e'$.

Applying Theorem \ref{thm:DF-2} with $\lambda=(p-1)/4$, we have the following theorem, which generalizes the latter part of Theorem 16.72(2) in \cite{ab}. Actually the former part of Theorem 16.72(2) can be also obtained by applying the following Theorem \ref{thm:DF-4} with $\lambda=(p-1)/2$.

\begin{theorem}\label{thm:cor}
There exists an $(\mathbb{F}_{p}\times \mathbb{F}_{q},\mathbb{F}_{p}\times \{0\},p,(p-1)/4)$-DF for any prime powers $p$ and $q$ with $p\equiv q\equiv 1 \pmod{4}$ and $q\geq p$.
\end{theorem}

\begin{theorem}\label{eg:17_2}
There exists an $(\mathbb{F}_{17}\times\mathbb{F}_q,\mathbb{F}_{17}\times\{0\},17,2)$-DF for all primes $q\equiv 1 \pmod{8}$.
\end{theorem}

\Proof Applying Theorem \ref{thm:DF-2} with $\lambda=2$ and $p=17$, we get an $(\mathbb{F}_{17}\times\mathbb{F}_q,\mathbb{F}_{17}\times\{0\},17,2)$-DF for all primes $q\equiv 1 \pmod{8}$ and $q>Q(2,13)=2.03024\times 10^9$. For primes $q\equiv 1 \pmod{8}$ and $q\leq Q(2,13)$, to apply Lemma \ref{lem:DF-2}, we have checked by computer search that there exist $y_1, y_2, y_3, y_4$ in $\mathbb{F}_q^*$ such that each $D_h$, $h\in\mathbb{F}_{17}$, is a $2$-transversal for the cosets of $C_0^{2,q}$ in $\mathbb{F}_q^*$. The interested reader may get a copy of these data from the authors. \qed

\subsection{The use of Lemma \ref{lem:DF-1} with $e=(q-1)/2$}

\begin{lemma}\label{lem:DF-4}
Let $p\equiv 1 \pmod{2}$ be a prime power and $\lambda$ be a divisor of $(p-1)/2$. Write $d=(p-1)/2\lambda$. Let $q$ be a prime power satisfying $\lambda(q-1)\equiv 0\pmod{p-1}$. Let $\delta$ be a generator of $C_0^{2,p}$. Take the first type Paley $(\mathbb{F}_{p},p,p-1)$-SDF whose unique base block $(f_{0},f_1,\ldots,f_{p-1})=$
\begin{eqnarray}\label{2st Paley}
(0,\delta,\delta,\delta^2,\delta^2,\ldots,\delta^{\frac{p-1}{2}},
\delta^{\frac{p-1}{2}}).
\end{eqnarray}
Suppose that one can choose an appropriate multiset $(\phi_{0},\phi_1,\ldots,\phi_{p-1})=$
\begin{eqnarray}\label{3st Paley-y}
(0,y_1,-y_1,y_2,-y_2,\ldots,y_{\frac{p-1}{2}},-y_{\frac{p-1}{2}})
\end{eqnarray}
such that $\{y_1,y_2,\ldots,y_{(p-1)/2}\}\subseteq {\mathbb F}_{q}^*$ and for each $h\in \mathbb{F}_{p}$,
\begin{eqnarray*}\label{3st Paley-dh}
[\phi_{a}-\phi_{b}:f_{a}-f_{b}=h,(a,b)\in I_{p}\times I_{p},a\neq b]=\{1,-1\}\cdot D_h,
\end{eqnarray*}
where $D_h$ is a $\lambda$-transversal for the cosets of $C_0^{d,q}$ in $\mathbb{F}_{q}^*$. Let $S$ be a representative system for the cosets of $\{1,-1\}$ in $C_0^{d,q}$. Let
$
B=\{(f_0,\phi_0),(f_1,\phi_1),\ldots,(f_{p-1},\phi_{p-1})\}.
$
Then
$$\mathfrak{F}=[B\cdot\{(1,s)\}:s\in S]$$
forms an $(\mathbb{F}_{p}\times \mathbb{F}_{q},\mathbb{F}_{p}\times \{0\},p,\lambda)$-DF. \end{lemma}

\Proof Apply Lemma \ref{lem:DF-1} with $n=1$, $k=p$, $d=(p-1)/2\lambda$ and $e=(q-1)/2$ to obtain the required $(\mathbb{F}_{p}\times \mathbb{F}_{q},\mathbb{F}_{p}\times \{0\},p,\lambda)$-DF. Note that $\{1,-1\}=C_0^{(q-1)/2,q}$. \qed

Similar arguments to those in Lemmas \ref{lem:DF-2-Dh} and \ref{lem:DF-2-Dh1}, we have

\begin{lemma}\label{lem:DF-4-Dh}
Follow the notation in Lemma $\ref{lem:DF-4}$.
\begin{itemize}
\item[$(1)$] W.l.o.g., $D_0=2\cdot\{y_1,y_2,\ldots,y_{(p-1)/2}\}$.
\item[$(2)$]
For each $h\in \mathbb{F}_{p}^*$, let
$$T_h=[\phi_{a}-\phi_{b}:f_{a}-f_{b}=h,(a,b)\in I_{p}\times I_{p},a\neq b].$$
Then $T_h=\{1,-1\}\cdot D_h$ for some $D_h\subset \mathbb{F}_{q}$ and the size of $D_h$ is $(p-1)/2$. Furthermore, $D_h=D_{-h}$ and w.l.o.g., $D_h$ consists of elements having types $(I)$ $y_i-y_j$, $y_i+y_j$, and $(II)$ $y_i$.
\end{itemize}\end{lemma}

\begin{theorem}\label{thm:DF-4}
Let $\lambda$ be any divisor of $(p-1)/2$ and $d=(p-1)/2\lambda$. There exists an $(\mathbb{F}_{p}\times \mathbb{F}_{q},\mathbb{F}_{p}\times \{0\},p,\lambda)$-DF for any prime powers $p$ and $q$ with $p\equiv 1 \pmod{2}$, $\lambda(q-1)\equiv 0\pmod{p-1}$ and $q>Q(d,p-2)$.
\end{theorem}

\Proof Combine the results of Lemmas \ref{lem:DF-4} and \ref{lem:DF-4-Dh}. Then apply Theorem \ref{thm:cyclot bound} as we have done in the proof of Theorem \ref{thm:DF-2} to complete the proof. \qed

We remark that when $p\in\{3,5\}$, by Theorems \ref{thm:DF-2} and \ref{thm:DF-4} one may refind the known results on $(pq,p,p,1)$-DFs over $\mathbb{F}_{p} \times \mathbb{F}_{q}$ listed in Theorem 16.71 in \cite{ab}.

\begin{lemma}\label{lem:DF-3}
Let $p\equiv 3 \pmod{4}$ be a prime power and $\lambda$ be a divisor of $(p+1)/2$. Write $d=(p+1)/2\lambda$. Let $q$ be a prime power satisfying $\lambda(q-1)\equiv 0\pmod{p+1}$. Let $\delta$ be a generator of $C_0^{2,p}$. Take the second type Paley $(\mathbb{F}_{p},p+1,p+1)$-SDF whose unique base block $(f_{0},f_1,\ldots,f_{p})=$
\begin{eqnarray}\label{2st Paley}
(0,0,\delta,\delta,\delta^2,\delta^2,\ldots,\delta^{\frac{p-1}{2}},
\delta^{\frac{p-1}{2}}).
\end{eqnarray}
Suppose that one can choose an appropriate multiset $(\phi_{0},\phi_1,\ldots,\phi_{p})=$
\begin{eqnarray}\label{2st Paley-y}
(y,-y,y_1,-y_1,y_2,-y_2,\ldots,y_{\frac{p-1}{2}},-y_{\frac{p-1}{2}})
\end{eqnarray}
such that $\{y,y_1,y_2,\ldots,y_{(p-1)/2}\}\subseteq {\mathbb F}_{q}^*$ and for each $h\in \mathbb{F}_{p}$,
\begin{eqnarray*}\label{2st Paley-dh}
[\phi_{a}-\phi_{b}:f_{a}-f_{b}=h,(a,b)\in I_{p+1}\times I_{p+1},a\neq b]=\{1,-1\}\cdot D_h,
\end{eqnarray*}
where $D_h$ is a $\lambda$-transversal for the cosets of $C_0^{d,q}$ in $\mathbb{F}_{q}^*$. Let $S$ be a representative system for the cosets of $\{1,-1\}$ in $C_0^{d,q}$. Let
$
B=\{(f_0,\phi_0),(f_1,\phi_1),\ldots,(f_{p},\phi_{p})\}.
$
Then
$$\mathfrak{F}=[B\cdot\{(1,s)\}:s\in S]$$
forms an $(\mathbb{F}_{p}\times \mathbb{F}_{q},\mathbb{F}_{p}\times \{0\},p+1,\lambda)$-DF. \end{lemma}

\Proof Apply Lemma \ref{lem:DF-1} with $n=1$, $k=p+1$, $d=(p+1)/2\lambda$ and $e=(q-1)/2$ to obtain the required $(\mathbb{F}_{p}\times \mathbb{F}_{q},\mathbb{F}_{p}\times \{0\},p,\lambda)$-DF. \qed

Similar arguments to those in Lemmas \ref{lem:DF-2-Dh} and \ref{lem:DF-2-Dh1}, we have

\begin{lemma}\label{lem:DF-3-Dh}
Follow the notation in Lemma $\ref{lem:DF-3}$.
\begin{itemize}
\item[$(1)$] W.l.o.g., $D_0=2\cdot\{y,y_1,y_2,\ldots,y_{(p-1)/2}\}$.
\item[$(2)$]
For each $h\in \mathbb{F}_{p}^*$, let
$$T_h=[\phi_{a}-\phi_{b}:f_{a}-f_{b}=h,(a,b)\in I_{p+1}\times I_{p+1},a\neq b].$$
Then $T_h=\{1,-1\}\cdot D_h$ for some $D_h\subset \mathbb{F}_{q}$ and the size of $D_h$ is $(p+1)/2$. Furthermore, $D_h=D_{-h}$ and w.l.o.g., $D_h$ consists of elements having types $(I)$ $y_i-y_j$, $y_i+y_j$, and $(II)$ $y_i-y$, $y_i+y$.
\end{itemize}\end{lemma}

\begin{theorem}\label{thm:DF-3}
Let $\lambda$ be any divisor of $(p+1)/2$ and $d=(p+1)/2\lambda$. There exists an $(\mathbb{F}_{p}\times \mathbb{F}_{q},\mathbb{F}_{p}\times \{0\},p+1,\lambda)$-DF for any prime powers $p$ and $q$ with $p\equiv 3 \pmod{4}$, $\lambda(q-1)\equiv 0\pmod{p+1}$ and $q>Q(d,p)$.
\end{theorem}

\Proof Combine the results of Lemmas \ref{lem:DF-3} and \ref{lem:DF-3-Dh}. Then apply Theorem \ref{thm:cyclot bound} as we have done in the proof of Theorem \ref{thm:DF-2} to complete the proof. \qed

\subsection{New $2$-designs}

Start from the relative difference families in Theorems \ref{thm:DF-2} and \ref{thm:DF-4}. Then apply Proposition \ref{prop:RDF-BIBD-1}(1) with a trivial $2$-$(p,p,\lambda)$ design. We obtain the following two theorems.

\begin{theorem}\label{thm:BIBD-1}
Let $\lambda$ be any divisor of $(p-1)/4$ and $d=(p-1)/4\lambda$.
\begin{itemize}
\item[$(1)$] There exists a $2$-$(pq,p,\lambda)$ design for any prime powers $p$ and $q$ with $p\equiv 1,5 \pmod{12}$, $\lambda(q-1)\equiv 0\pmod{p-1}$ and $q>Q(d,p-4)$.
\item[$(2)$] Let $p$ be a power of $9$, and $\xi$ be a primitive $4$th root of unity in $\mathbb{F}_{q}$. If $\lambda>1$ or $1-\xi\not\in C_{0}^{d,q}$, then there exists a $2$-$(pq,p,\lambda)$ design for any prime power $q$ with $\lambda(q-1)\equiv 0\pmod{p-1}$ and $q>Q(d,p-4)$.
\end{itemize}
\end{theorem}

\begin{theorem}\label{thm:BIBD-2}
Let $\lambda$ be any divisor of $(p-1)/2$ and $d=(p-1)/2\lambda$. There exists a $2$-$(pq,p,\lambda)$ design for any prime powers $p$ and $q$ with $p\equiv 1 \pmod{2}$, $\lambda(q-1)\equiv 0\pmod{p-1}$ and $q>Q(d,p-2)$.
\end{theorem}

Start from the relative difference families in Theorem \ref{thm:DF-3}. Then apply Proposition \ref{prop:RDF-BIBD-1}(2) with a trivial $2$-$(p+1,p+1,\lambda)$ design. We obtain the following theorem.

\begin{theorem}\label{thm:BIBD-3}
Let $\lambda$ be any divisor of $(p+1)/2$ and $d=(p+1)/2\lambda$. There exists a $2$-$(pq+1,p+1,\lambda)$ design for any prime powers $p$ and $q$ with $p\equiv 3 \pmod{4}$, $\lambda(q-1)\equiv 0\pmod{p+1}$ and $q>Q(d,p)$.
\end{theorem}

\begin{remark}\label{remark:cyclic}
If $p$ and $q$ are both prime and $p\neq q$, then since $\mathbb{F}_{p}\times\mathbb{F}_{q}$ is isomorphic to $\mathbb{Z}_{pq}$, all $2$-$(pq,p,\lambda)$ designs from Theorems $\ref{thm:BIBD-1}$ and $\ref{thm:BIBD-2}$ are cyclic, while all $2$-$(pq+1,p+1,\lambda)$ designs from Theorem $\ref{thm:BIBD-3}$ are $1$-rotational.
\end{remark}

Combine the results of Theorems \ref{eg:13_1}, \ref{eg:17_1}, \ref{thm:cor} and \ref{eg:17_2}, and apply Proposition \ref{prop:RDF-BIBD-1}(1). We have

\begin{theorem}\label{thm:13_17}
\begin{itemize}
\item[$(1)$] There exists a $2$-$(13q,13,1)$ design for any prime $q\equiv 1 \pmod{12}$ with the possible exceptions of $q\in E_{13}=\{37, 61, 73, 97, 109, 181, 313, 337, 349, 373, $ $409, 421, 541, 577, 829, 853, 1129, 1741, 2473\}.$
\item[$(2)$] There exists a $2$-$(17q,17,1)$ design for any prime $q\equiv 1 \pmod{16}$ and $q>Q(4,13)$ $=3.44807\times 10^{17}$, or $q\in S_{17}\cup\{p: p\ {\rm is\ a\ prime}, p\equiv1\pmod{16}, 6673\leq p\leq 9857\}$ $($see Theorem $\ref{eg:17_1}$ for details of $S_{17})$.
\item[$(3)$] There exists a $2$-$(17q,17,2)$ design for any prime $q\equiv 1 \pmod{8}$.
\item[$(4)$] There exists a $2$-$(pq,p,(p-1)/4)$ design for any prime powers $p$ and $q$ with $p\equiv q\equiv 1 \pmod{4}$ and $q\geq p$.
\end{itemize}
\end{theorem}

We remark that M. Buratti discussed constructions for $2$-$(13q,13,\lambda)$ designs and $2$-$(17q,17,\lambda)$ designs in \cite{b97jcta,b99dcc}. His results rely on cyclotomic conditions of some specific elements. For example, to construct $2$-$(17q,17,2)$ designs, 2 is required to be not a 4th power in $\mathbb{F}_q$ \cite{b99dcc}. Theorem \ref{thm:13_17}(4) generalizes Theorem 4.7 in \cite{b99}.

\begin{remark}\label{remark:13 17}
If $q\neq 13$ in Theorem $\ref{thm:13_17}(1)$, or $q\neq 17$ in Theorem $\ref{thm:13_17}(2)$ and $(3)$, or $q\neq p$ in Theorem $\ref{thm:13_17}(4)$, then all $2$-designs from Theorem $\ref{thm:13_17}$ are cyclic.
\end{remark}

\section{Concluding remarks}

By a careful application of cyclotomic conditions attached to strong difference families, this paper improves the lower bound on the asymptotic existence results of $(\mathbb{F}_{p}\times \mathbb{F}_{q},\mathbb{F}_{p}\times \{0\},k,\lambda)$-DFs for $k\in\{p,p+1\}$, and presents seven new $2$-designs.

Future directions are two-fold. One is to systematically analyze possible patterns of SDFs that help to produce new 2-designs from the point of view of asymptotic existence or concrete designs. We point out that M. Buratti et al. \cite{bcw,brz} essentially made use of a Paley SDF of the first type, called a strong difference map from a graph-theoretical perspective, to investigate the constructions for $i$-perfect cycle decompositions.

The other direction is to generalize the techniques used in this paper to construct other kinds of difference families, such as resolvable difference families (cf. \cite{b97}) and partitioned difference families (cf. \cite{byw}), which can be used to construct frequency hopping sequences (cf. \cite{gfm}) and constant composition codes (cf. \cite{yst}), etc.

\subsection*{Acknowledgements}
The authors wishes to express their sincere appreciations to all those who made suggestions for improvements to this paper. Particularly thanks go to Professor M. Buratti who critically read the paper and made numerous helpful suggestions including pointing out the importance of Lemma \ref{lem:DF-1} and simplifying the proof of Theorem \ref{thm:DF-2}.

\end{document}